\documentstyle[12pt]{article}

\font\twlgot =eufm10 scaled \magstep1
\font\egtgot =eufm8
\font\sevgot =eufm7
\font\twlmsb =msbm10 scaled \magstep1
\font\egtmsb =msbm8
\font\sevmsb =msbm7

\newfam\gotfam

\textfont\gotfam\twlgot
\scriptfont\gotfam\egtgot
\scriptscriptfont\gotfam\sevgot

\newfam\msbfam
\textfont\msbfam\twlmsb
\scriptfont\msbfam\egtmsb
\scriptscriptfont\msbfam\sevmsb
\def\Bbb{\protect\pBbb}
\def\pBbb{\relax\ifmmode\expandafter\Bb\else\typeout{You cann't use
Bbb in text mode}\fi}
\def\Bb #1{{\fam\msbfam\relax#1}}

\def\thebibliography#1{\section*{References}\list
  {[\arabic{enumi}]}{\settowidth\labelwidth{#1}\leftmargin\labelwidth
    \advance\leftmargin\labelsep
    \usecounter{enumi}}
    \def\newblock{\hskip .11em plus .33em minus .07em}
    \sloppy\clubpenalty4000\widowpenalty4000
    \sfcode`\.=1000\relax}

\def\op#1{\mathop{\fam0 #1}\limits}

\newcommand{\beq}{\begin{equation}}
\newcommand{\eeq}{\end{equation}}
\newcommand{\ben}{\begin{eqnarray}}
\newcommand{\een}{\end{eqnarray}}
\newcommand{\be}{\begin{eqnarray*}}
\newcommand{\ee}{\end{eqnarray*}}
\newcommand{\bea}{\begin{eqalph}}
\newcommand{\eea}{\end{eqalph}}

\newcommand{\cH}{{\cal H}}
\newcommand{\cF}{{\cal F}}

\newcommand{\cG}{{\cal G}}

\newcommand{\al}{\alpha}

\newcommand{\la}{\lambda}

\newcommand{\Om}{\Omega}
\newcommand{\m}{\mu}

\newcommand{\g}{\gamma}

\newcommand{\vt}{\vartheta}
\newcommand{\vf}{\varphi}
\newcommand{\up}{\upsilon}

\newcommand{\di}{{\rm dim\,}}

\newcommand{\w}{\wedge}

\newcommand{\dr}{\partial}
\newcommand{\ar}{\op\longrightarrow}

\newcommand{\ve}{\varepsilon}

\newcounter{eqalph}
\newcounter{equationa}
\newcounter{remark}
\newcounter{example}
\newcounter{theorem}
\newcounter{proposition}
\newcounter{lemma}
\newcounter{corollary}
\newcounter{definition}
\def\theremark{\arabic{remark}}

\def\thedefinition{\arabic{definition}}

\newenvironment{theo}{\refstepcounter{definition}
\bigskip\noindent{\bf Theorem \thedefinition.} \it}{\medskip}

\newenvironment{lem}{\refstepcounter{definition}
\bigskip\noindent{\bf Lemma \thedefinition.}\it}{\medskip}

\newenvironment{eqalph}{\stepcounter{equation}
\setcounter{equationa}{\value{equation}} \setcounter{equation}{0}

\begin{eqnarray}}{\end{eqnarray}
\setcounter{equation}{\value{equationa}}}

\newcommand{\mar}[1]{}

\begin{document}
\hbox{}

{\parindent=0pt

{\large \bf Noncommutative integrability on noncompact invariant
manifolds}
\bigskip

{\bf E.Fiorani}$^1$, {\bf G.Sardanashvily}$^2$
\bigskip

\begin{small}

$^1$ Department of Mathematics and Informatics, University of
Camerino, 62032 Camerino (MC), Italy

\medskip

$^2$ Department of Theoretical Physics, Moscow State University,
117234 Moscow, Russia

\bigskip

{\bf Abstract}

The Mishchenko--Fomenko theorem on noncommutative completely integrable
Hamiltonian systems on a symplectic manifold is extended to the case of
noncompact invariant submanifolds.

\end{small}
\bigskip


 }

\section{Introduction}

We are concerned with the classical theorems on abelian and
noncommutative integrability of Hamiltonian systems on a
symplectic manifold. These are the Liouville--Arnold theorem on
abelian completely integrable systems (henceforth CIS)
\cite{arn,arn1,laz},  the Poincar\'e--Lyapounov--Nekhoroshev
theorem on abelian partially integrable systems
\cite{gaeta,nekh,nekh94}, and the Mishchenko--Fomenko one on
nocommutative CISs \cite{fasso05,karas,mishc}. These theorems
state the existence of (generalized) action-angle coordinates
around a compact invariant submanifold, which is a torus. However,
there is a topological obstruction to the existence of global
action-angle coordinates \cite{daz,dust}. The Liouville--Arnold
and Nekhoroshev theorems have been extended to noncompact
invariant submanifolds, which are toroidal cylinders
\cite{fior,fior1,jmp03,vin}.  In particular, this is the case of
time-dependent CISs \cite{jpa02,book05}. Any time-dependent CIS of
$m$ degrees of freedom can be represented as an autonomous one of
$m+1$ degrees of freedom on a homogeneous momentum phase space,
where time is a generalized angle coordinate. Therefore, we
further consider autonomous CISs.

Our goal here is the following generalization of the
Mishchenko--Fomenko theorem to noncommutative CISs whose invariant
submanifolds need not be compact.

\begin{theo} \label{nc0} \mar{nc0} Let $(Z,\Om)$ be a
connected symplectic $2n$-dimensional real smooth manifold and
$(C^\infty(Z),
\{,\})$ the Poisson algebra of smooth real functions on $Z$. Let a subset
$H=(H_1,\ldots,H_k)$,
$n\leq k<2n$, of $C^\infty(Z)$ obey the
following conditions.

(i) The Hamiltonian vector fields $\vt_i$ of functions $H_i$ are complete.

(ii) The map $H:Z\to \Bbb R^k$ is a submersion with connected and
mutually diffeomorphic fibers, i.e.,
\mar{nc4}\beq
H:Z\to N=H(Z) \label{nc4}
\eeq
is a fibered manifold over a connected open subset $N\subset\Bbb R^k$.

(iii)  There exist real smooth functions $s_{ij}: N\to \Bbb R$ such that
\mar{nc1}\beq
\{H_i,H_j\}= s_{ij}\circ H, \qquad i,j=1,\ldots, k. \label{nc1}
\eeq

(iv) The matrix function with the entries $s_{ij}$ (\ref{nc1}) is of
constant corank
$m=2n-k$ at all points of $N$.

\noindent
Then the following hold.

(I) The fibers of $H$ (\ref{nc4}) are diffeomorphic to a toroidal
cylinder
\mar{nc2}\beq
\Bbb R^{m-r}\times T^r. \label{nc2}
\eeq

(II) Given a fiber $M$ of $H$ (\ref{nc4}), there exists an open
saturated neighbourhood $U_M$ of it (i.e., a fiber through a point
of $U_M$ belongs to $U_M$), which is a trivial principal bundle
with the structure group (\ref{nc2}).

(III) Given standard coordinates $(y^\la)$ on
the toroidal cylinder (\ref{nc2}), the neighbourhood $U_M$ is provided
with bundle coordinates
$(J_\la,p_A,q^A,y^\la)$, called the generalized action-angle coordinates,
which are the Darboux coordinates of the symplectic form
$\Om$ on $U_M$, i.e.,
\mar{nc3}\beq
\Om= dJ_\la\w dy^\la + dp_A\w dq^A. \label{nc3}
\eeq
\end{theo}

In Hamiltonian mechanics, one can think of functions $H_i$
in Theorem \ref{nc0} as being integrals of motion of a CIS. Their
level surfaces (fibers of $H$) are invariant submanifolds of a CIS.

\section{Abelian completely and partially integrable systems}

The proof of Theorem 1 is based on the fact that an invariant
submanifold of a noncommutative CIS is a maximal integral manifold
of some abelian partially integrable Hamiltonian system
\cite{fasso05}.

If $k=n$, Theorem \ref{nc0} provides the above mentioned extension
of the Liouville--Arnold theorem to abelian CISs whose invariant
submanifolds are noncompact (\cite{vin}, Theorem 6.1; \cite{fior},
Theorem 1). Note that the proof of Theorem 6.1 \cite{vin} differs
from that of Theorem 1 \cite{fior}. It is based on Lemma 6.4. The
statement of its Corollary 6.3 is equivalent to the assumption of
Lemmas 6.1 -- 6.4 that an imbedded invariant submanifold
$N_x\subset K$ admits a Lagrangian transversal submanifold
$W\subset K$ through $x$. Apparently, one can avoid the
construction of $T^*(W)$ from the proof and, instead of the map
$\g$, consider the map
\be
W\times \Bbb R^n\ar (W\times \Bbb R^n)/\Bbb Z^{m(x)}\ar
\al(W\times \Bbb R^n).
\ee
One also need not appeal to the concept of many-valued functions
$\vf_i$, but can show that the fibered manifold $\al(W\times \Bbb
R^n)\to W$ is a fiber bundle. This is always true if its fibers
are tori and, if $m(x)<n$, follows from the fact that sections
$l_i$ of $W\times \Bbb R^n\to W$ introduced in Lemma 6.2 are
smooth.

The condition (ii) of Theorem \ref{nc0} implies that the functions
$\{H_\la\}$ are independent on $Z$, i.e., the $n$-form $\op\w^n
dH_\la$ nowhere vanishes. Accordingly, the Hamiltonian vector
fields $\vt_\la$ of these functions are independent on $Z$, i.e.,
the multivector field $\op\w^n\vt_\la$ nowhere vanishes. If $k=n$,
these vector fields are mutually commutative, and they span a
regular involutive $n$-dimensional distribution on $Z$ whose
maximal integral manifolds are exactly fibers of the fibered
manifold (\ref{nc4}). Thus, every fiber of $H$ (\ref{nc4}) admits
$n$ independent complete vector fields, i.e., it is a locally
affine manifold and, consequently, diffeomorphic to a toroidal
cylinder.

Considering an abelian CIS around some compact invariant
submanifold, we come to the Liouville--Arnold theorem (somebody
also calls it the Liouville--Mineur--Arnold theorem \cite{zung}).
Instead of the conditions (i) and (ii) of Theorem 1, one can
suppose that integrals of motion $\{H_\la\}$ are independent
almost everywhere on a symplectic manifold $Z$, i.e., the set of
points where the exterior form $\op\w^n dH_\la$ (or, equivalently,
the multivector field $\op\w^n\vt_\la$) vanishes is nowhere dense.
In this case, connected components of level surfaces of functions
$\{H_\la\}$ form a singular Stefan foliation $\cF$ of $Z$ whose
leaves are both the maximal integral manifolds of the singular
involutive distribution spanned by the vector fields $\vt_\la$ and
the orbits of the pseudogroup $G$ of local diffeomorphisms of $Z$
generated by flows of these vector fields \cite{stef,susm}. Let
$M$ be a leaf of $\cF$ through a regular point $z\in Z$ where
$\op\w^n\vt_\la\neq 0$. It is regular everywhere because the group
$G$  preserves $\op\w^n\vt_\la$. If $M$ is compact and connected,
there exists its saturated open neighbourhood $U_M$ such that the
map $H$ restricted to $U_M$ satisfies the condition (ii) of
Theorem 1, i.e., the foliation $\cF$ of $U_M$ is a fibered
manifold in tori $T^n$. Since its fibers are compact, $U_M$ is a
bundle \cite{meig}. Hence, it contains a saturated open
neighbourhood of $M$, say again $U_M$, which is a trivial
principal bundle with the structure group $T^n$. Providing $U_M$
with the Darboux (action-angle) coordinates $(J_\la,\al^\la)$, one
uses the fact that there are no linear functions on a torus $T^n$.

The Poincar\'e--Lyapounov--Nekhoroshev theorem generalizes the
Liouville--Arnold one to partially integrable systems
characterized by $k<n$ independent integrals of motion $H_\la$ in
involution. In this case, one deals with $k$-dimensional maximal
integral manifolds of the distribution spanned by Hamiltonian
vector fields $\vt_\la$ of integrals of motion $H_\la$. The
Poincar\'e--Lyapounov--Nekhoroshev theorem imposes a sufficient
condition which Hamiltonian vector fields $\vt_\la$ must satisfy
in order that their compact maximal integral manifold $M$ admits
an open neighbourhood fibered in tori \cite{gaeta,gaeta03}. Such a
condition has been also investigated in the  case of
noncommutative vector fields depending on parameters
\cite{gaeta06}. Extending the Poincar\'e--Lyapounov--Nekhoroshev
theorem to the case of noncompact integral submanifolds, we in
fact assumed from the beginning that these submanifolds form a
fibration \cite{fior,jmp03,book05}. In a more general setting, we
have studied the property of a given dynamical system to be
Hamiltonian relative to different Poisson structures
\cite{bols,jmp03,spar}. As is well known, any integrable
Hamiltonian system is Hamiltonian relative to different symplectic
and Poisson structures, whose variety has been analyzed from
different viewpoints
\cite{bog96,bog98,brouz,fasso98,fasso02,magri,smirn}. One of the
reasons is that bi-Hamiltonian systems have a large supply of
integrals of motion.   Here, we refer to our following result on
partially integrable systems on a symplectic manifold
(\cite{jmp03}, Theorem 6).

\begin{theo} \label{nc6} \mar{nc6}
Given a $2n$-dimensional symplectic manifold $(Z,\Om)$, let
$\{H_1,\ldots,H_m\}$, $m\leq n$, be smooth real functions on $Z$ in
involution which satisfy the following conditions.

(i) The functions $H_\la$ are everywhere independent.

(ii) Their Hamiltonian vector fields $\vt_\la$ are complete.

(iii) These vector fields span a regular distribution whose maximal
integral manifolds form a fibration $\cF$ of $Z$ with diffeomorphic
fibers.

\noindent
Then the following hold.

(I) All fibers of $\cF$ are diffeomorphic to a toroidal cylinder
(\ref{nc2}).

(II) There is an open saturated neighbourhood $U_M$ of every fiber
$M$ of $\cF$ which is a trivial principal bundle with the
structure group (\ref{nc2}).

(III) Given standard coordinates $(y^\la)$ on the toroidal
cylinder (\ref{nc2}), the neighbourhood $U_M$ is endowed with the
bundle coordinates $(J_\la,p_A,q^A,y^\la)$ such that the symplectic
form $\Om$ is brought into the form (\ref{nc3}).
\end{theo}

Theorem \ref{nc6} provides the above mentioned generalization of
the Poincar\'e--Lyapounov--Nekhoroshev theorem to the case of
noncompact invariant submanifolds. A geometric aspect of this
generalization is the following.  Any fibered manifold whose
fibers are diffeomorphic either to $\Bbb R^r$ or a compact
connected manifold $K$ (e.g., a torus) is a fiber bundle
\cite{meig}. However, a fibered manifold whose fibers are
diffeomorphic to a product $\Bbb R^r\times K$ (e.g., a toroidal
cylinder (\ref{nc2})) need not be a fiber bundle (see \cite{book},
Example 1.2.2).

\section{The proof of Theorem 1}

Theorem \ref{nc6} is the final step of the proof of Theorem
\ref{nc0}. The condition (iv) of Theorem \ref{nc0} implies that an
$m$-dimensional invariant submanifold of a noncommutative CIS is a
maximal integral manifold of some abelian partially integrable
Hamiltonian system obeying the conditions of Theorem \ref{nc6}.
The proof of this fact is based on the following two assertions
\cite{fasso05,libe}.

\begin{lem} \label{nc7} \mar{nc7} Given a symplectic manifold $(Z,\Om)$,
let $H:Z\to N$ be a fibered manifold such
that, for any two functions $f$, $f'$ constant on fibers of $H$, their
Poisson bracket $\{f,f'\}$ is so. Then $N$ is provided with an
unique coinduced Poisson structure $\{,\}_N$ such that $H$ is a Poisson
morphism.
\end{lem}

Since any function constant on fibers of $H$ is a pull-back of some
function on $N$, the condition of Lemma \ref{nc7} is satisfied due to
item (iii) of Theorem \ref{nc0}. Thus, the base $N$ of the fibration
(\ref{nc4}) is endowed with a coinduced Poisson structure.

\begin{lem} \label{nc8} \mar{nc8} Given a fibered manifold $H:Z\to N$,
the following conditions are equivalent:

(i) the rank of the coinduced Poisson structure $\{,\}_N$ on $N$ equals
$2\di N-\di Z$,

(ii) the fibers of $H$ are isotropic,

(iii) the fibers of $H$ are  maximal integral manifolds of the involutive
distribution spanned by the Hamiltonian vector fields of the pull-back
$H^*C$ of Casimir functions $C$ of the Poisson algebra on $N$.
\end{lem}

It is readily observed that the condition (i) of Lemma \ref{nc8}
is satisfied due to the assumption (iv) of Theorem \ref{nc0}. It
follows that every fiber $M$ of the fibration (\ref{nc4}) is a
maximal integral manifold of the involutive distribution spanned
by the Hamiltonian vector fields $\up_\la$ of the pull-back
$H^*C_\la$ of $m$ independent Casimir functions $\{C_1,\ldots,
C_m\}$ on an open neighbourhood $N_M$ of the point $H(M)$. Let us
put $U_M=H^{-1}(N_M)$. It is an open saturated neighbourhood of
$M$. Since
\mar{j21}\beq
H^*C_\la(z)= (C_\la\circ
H)(z)= C_\la(H_i(z)),\qquad z\in U_M, \label{j21}
\eeq
the Hamiltonian vector fields $\up_\la$ on $M$ are linear combinations
of Hamiltonian vector fields $\vt_i$ of the functions $H_i$ and,
therefore, they are complete on $M$. Similarly, they are complete on any
fiber of $U_M$ and, consequently, on $U_M$. Thus, the conditions of
Theorem \ref{nc6} hold on $U_M$. This completes the proof of Theorem
\ref{nc0}.

The proof of Theorem \ref{nc0} gives something more. Let $\{H_i\}$
be integrals of motion of a Hamiltonian $\cH$. Since
$(J_\la,p_A,q^A)$ are coordinates on $N$, they are also integrals
of motion of $\cH$. Therefore, the Hamiltonian $\cH$ depends only
on the action coordinates $J_\la$, and the equation of motion read
\be
\dot y^\la=\frac{\dr\cH}{\dr J_\la}, \quad J_\la={\rm const.},
\quad q^A={\rm const.}, \quad p_A={\rm const.}
\ee
Though the integrals of
motion
$H_i$ are smooth functions of coordinates
$(J_\la,q^A,p_A)$, the Casimir functions
\be
C_\la(H_i(J_\m,q^A,p_A))=
C_\la(J_\m)
\ee
depend only on the  action coordinates $J_\la$. Moreover, a
Hamiltonian
\be
\cH(J_\m)=\cH(C_\la(H_i(J_\m,q^A,p_A))
\ee
is expressed in integrals of
motion $H_i$ through the Casimir functions (\ref{j21}).

Let us note that, under the assumptions of the
Mishchenko--Fomenko theorem, a noncommutative CIS is also
integrable in the abelian sense. Namely, it admits $n$
independent integrals of motion in involution \cite{bols03}. Under the
conditions of Theorem \ref{nc0}, such
integrals of motion in involution exist, too. All of them are
the pull-back of functions on $N$. However, one must justify that they
obey the condition (iii) of Theorem
\ref{nc6} in order to characterize an abelian CIS.

\section{Example}

The original Mishchenko--Fomenko theorem is restricted to CISs
whose integrals of motion $\{H_1,\ldots,H_k\}$ form a $k$-dimensional real
Lie algebra $\cG$ of rank $m$ with the commutation relations
\be
\{H_i,H_j\}= c_{ij}^h H_h, \qquad c_{ij}^h={\rm const.}
\ee
In this case, nonvanishing complete Hamiltonian vector fields
$\vt_i$ of $H_i$ define a free Hamiltonian action on $Z$ of some
connected Lie group $G$ whose Lie algebra is isomorphic to $\cG$.
Orbits of $G$ coincide with $k$-dimensional maximal integral
manifolds of the regular distribution on $Z$ spanned by
Hamiltonian vector fields $\vt_i$ \cite{susm}. Furthermore, one
can treat $H$ (\ref{nc4}) as an equivariant momentum mapping of
$Z$ to the Lie coalgebra $\cG^*$, provided with the coordinates
$x_i(H(z))=H_i(z)$, $z\in Z$, \cite{book05,guil}. In this case,
the coinduced Poisson structure $\{,\}_N$ in Lemma \ref{nc7}
coincides with the canonical Lie--Poisson structure on $\cG^*$
given by the Poisson bivector field
\be
w=\frac12 c_{ij}^h x_h\dr^i\w\dr^j.
\ee
Recall that the coadjoint action of $\cG$ on $\cG^*$ reads
$\ve_i(x_j)=c_{ij}^hx_h$, where $\{\ve_i\}$ is a basis for $\cG$.
Casimir functions of the Lie--Poisson structure are exactly the
coadjoint invariant functions on $\cG^*$. They are constant on
orbits of the coadjoint action of $G$ on $\cG^*$ which coincide
with leaves of the symplectic foliation of $\cG^*$. Given a point
$z\in Z$ and the orbit $G_z$ of $G$ in $Z$ through $z$, the
fibration $H$ (\ref{nc4}) projects this orbit onto the orbit
$G_{H(z)}$ of the coadjoint action of $G$ in $\cG^*$ through
$H(z)$. Moreover, by virtue of item (iii), Lemma \ref{nc8}, the
inverse image $H^{-1}(G_{H(z)})$ of $G_{H(z)}$ coincides with the
orbit $G_z$. It follows that any orbit of $G$ in $Z$ is fibered in
invariant submanifolds.

The Mishchenko--Fomenko
theorem has been mainly applied to CISs whose integrals of motion form a
compact Lie algebra. Indeed, the group $G$ generated by flows of
the Hamiltonian vector fields is compact, and every orbit of $G$ in $Z$
is compact. Since a fibration of a compact manifold possesses compact
fibers, any invariant submanifold of such a noncommutative CIS is
compact.  Therefore, our Theorem
\ref{nc0} essentially extends a class of noncommutative CISs under
investigation.

For instance, a spherical top exemplifies a noncommutative CIS
whose integrals of motion make up the compact Lie algebra $so(3)$
with respect to some symplectic structure.

Let us consider
a CIS with the Lie algebra
$\cG=so(2,1)$ of integrals of motion $\{H_1,H_2,H_3\}$ on a
four-dimensional symplectic manifold
$(Z,\Om)$, namely,
\mar{j50}\beq
\{H_1,H_2\}=-H_3, \qquad \{H_2,H_3\}=H_1, \qquad \{H_3,H_1\}=H_2.
\label{j50}
\eeq
The rank of this algebra (the dimension of its Cartan subalgebra)
equals one.  Therefore, an invariant submanifold in Theorem
\ref{nc0} is $M=\Bbb R$, provided with a Cartesian coordinate $y$.
Let us consider its open saturated neighbourhood $U_M$ projected
via $H:U_M\to N$ onto a domain $N\subset \cG^*$ in the Lie
coalgebra $\cG^*$ centered at a point $H(M)\in \cG^*$ which
belongs to an orbit of the coadjoint action of maximal dimension
2. A domain $N$ is endowed with the coordinates $(x_1,x_2,x_3)$
such that integrals of motion  $\{H_1,H_2,H_3\}$ on $U_M=N\times
\Bbb R$, coordinated by $(x_1,x_2,x_3,y)$, read
\be
H_1=x_1, \qquad H_2=x_2, \qquad H_3=x_3.
\ee
As was mentioned above, the coinduced Poisson structure on $N$
is the Lie--Poisson structure
\mar{j51}\beq
w= x_2\dr^3\w\dr^1 - x_3\dr^1\w\dr^2 + x_1\dr^2\w\dr^3. \label{j51}
\eeq
Let us endow $N$ with different coordinates $(r,x_1,\g)$ given by
the equalities
\mar{j52}\beq
r=(x_1^2 + x_2^2 - x_3^2)^{1/2}, \quad x_2=(r^2-x_1^2)^{1/2}{\rm ch}\g,
\quad x_3=(r^2-x_1^2)^{1/2}{\rm sh}\g, \label{j52}
\eeq
where $r$ is a Casimir function on $\cG^*$. It is readily observed that
the coordinates (\ref{j52}) are the Darboux coordinates of the
Lie--Poisson structure (\ref{j51}), namely,
\mar{j53}\beq
w=\frac{\dr}{\dr \g}\w \frac{\dr}{\dr x_1}. \label{j53}
\eeq
Let $\vt_r$ be the Hamiltonian vector field of the Casimir function $r$
(\ref{j52}). This vector field is a combination
\be
\vt_r=\frac1{r}(x_1\vt_1 + x_2\vt_2 - x_3\vt_3)
\ee
of the Hamiltonian vector fields $\vt_i$ of integrals of motion $H_i$.
Its flows are invariant submanifolds. Let $y$ be a parameter
along the flows of this vector field, i.e.,
\be
\vt_r= \frac{\dr}{\dr y}.
\ee
Then the Poisson bivector associated to the symplectic form $\Om$
on $U_M$ is
\mar{j54}\beq
W= \frac{\dr}{\dr r}\w \frac{\dr}{\dr y} + \frac{\dr}{\dr \g}\w
\frac{\dr}{\dr x_1}. \label{j54}
\eeq
Accordingly, Hamiltonian vector fields of integrals of motion take
the form
\be
&& \vt_1= -\frac{\dr}{\dr \g}, \\
&& \vt_2= r(r^2-x_1^2)^{-1/2}{\rm ch}\g
\frac{\dr}{\dr y} + x_1 (r^2-x_1^2)^{-1/2}{\rm ch}\g
\frac{\dr}{\dr \g} + (r^2-x_1^2)^{1/2}{\rm sh}\g
\frac{\dr}{\dr x_1}, \\
&& \vt_3= r(r^2-x_1^2)^{-1/2}{\rm sh}\g
\frac{\dr}{\dr y} + x_1 (r^2-x_1^2)^{-1/2}{\rm sh}\g
\frac{\dr}{\dr \g} + (r^2-x_1^2)^{1/2}{\rm ch}\g
\frac{\dr}{\dr x_1}.
\ee

Thus, a symplectic annulus $(U_M, W)$ around an invariant
submanifold $M=\Bbb R$ is endowed with the generalized action-angle
coordinates
$(r,x_1,\g,y)$, and possesses the corresponding noncommutative CIS $\{r,
H_1,
\g\}$ with the commutation relations
\be
\{r, H_1\}=0, \qquad \{r, \g\}=0, \qquad \{H_1,\g\}=1.
\ee
This CIS
is related to the original one by the  transformations
\be
r=(H_1^2+ H_2^2 - H_3^2)^{1/2}, \qquad H_2=(r^2-H_1^2)^{1/2}{\rm ch}\g,
\qquad H_3=(r^2-H_1^2)^{1/2}{\rm ch}\g.
\ee
Its Hamiltonian is expressed only in the action
variable $r$.


\begin{thebibliography}{ddd}

\bibitem{arn} Arnold V and Avez A 1968 {\it Ergodic Problems in Classical Mechanics}
(Benjamin: NY)

\bibitem{arn1} Arnold V (Ed.) 1990 {\it Dynamical Systems III, IV}
(Springer: Berlin)

\bibitem{bog96} Bogoyavlenskij O 1996 Theory of tensor invariants of
integrable hamiltonian systems. I. Incompatible Poisson structures {\it
Commun. Math. Phys.} {\bf 180} 529-586

\bibitem{bog98} Bogoyavlenskij O 1998 Extended integrability and
bi-Hamiltonian systems {\it Commun. Math. Phys.}
{\bf 196} 19-51

\bibitem{bols} Bolsinov A and Borisov A 2002 Compatible Poisson brackets
on Lie algebras {\it Math. Notes} {\bf 72} 10-30

\bibitem{bols03} Bolsinov A and Jovanovi\'c B 2003 Noncommutative
integrability, moment map and geodesic flows {\it Ann. Global Anal.
Geom.} {\bf 23} 305-322

\bibitem{bols04} Bolsinov A and Jovanovi\'c B 2004 Complete involutive
algebras of functions on cotangent bundles of homogeneous spaces {\it
Math. Z} {\bf 246} 213-236

\bibitem{brouz} Brouzet R 1993  About the existence of recursion
operators for completely integrable
Hamiltonian systems near a Liouville torus
{\it J. Math. Phys.} {\bf 34} 1309-1313

\bibitem{cush} Cushman R and Bates L 1997 {\it Global Aspects of
Classical Integrable Systems} (Birkh\"auser: Basel)

\bibitem{daz} Dazord P and Delzant T 1987 Le probleme general des
variables actions-angles {\it J. Diff. Geom.} {\bf 26} 223-251

\bibitem{dust} Dustermaat J 1980 On global action-angle coordinates {\it
Commun. Pure Appl. Math.}
 {\bf 33} 687-706


\bibitem{fasso98} Fass\'o F and Ratiu T 1998 Compatibility of symplectic
structures adapted to noncommutatively integrable systems {\it J. Geom.
Phys.} {\bf 27} 199-220

\bibitem{fasso02} Fass\'o F and Giacobbe A 2002 Geometric structure
of "broadly integrable" vector
fields {\it J. Geom. Phys.} {\bf 44} 156-170

\bibitem{fasso05} Fass\'o F 2005 Superintegrable Hamiltonian systems:
geometry and applications {\it Acta Appl. Math.} {\bf 87} 93-121


\bibitem{fior} Fiorani E, Giachetta G and Sardanashvily G 2003 The
Liouville--Arnold--Nekhoroshev theorem for non-compact invariant
manifolds {\it J. Phys. A} {\bf 36} L101-L107

\bibitem{fior1} Fiorani E 2004 Completely and partially integrable systems
in the noncompact space {\it Int. J. Geom. Methods Mod. Phys.} {\bf 1}
167-183


\bibitem{gaeta} Gaeta G 2002 The Poincar\'e-Lyapunov-Nekhoroshev theorem
{\it Ann. Phys.} {\bf 297} 157-173

\bibitem{gaeta03} Gaeta G 2003 The Poincar\'e--Nekhoroshev map {\it J.
Nonlin. Math. Phys.} {\bf 10} 51-64

\bibitem{gaeta06} Gaeta G 2006 The Poincar\'e-Lyapunov-Nekhoroshev
theorem for involutory systems of vector fields {\it Ann. Phys.}
{\bf 321} 1277-1296

\bibitem{book} Giachetta G, Mangiarotti L and Sardanashvily G 1997 {\it
New Lagrangian and Hamiltonian Methods in Field Theory} (World
Scientific: Singapore)

\bibitem{jpa02} Giachetta G, Mangiarotti L. and Sardanashvily G 2002
Action-angle coordinates for time-dependent completely integrable Hamiltonian systems,
{\it J. Phys. A} {\bf 35} L439-L445

\bibitem{jmp03} Giachetta G, Mangiarotti L and Sardanashvily G 2003
Bi-Hamiltonian partially integrable systems {\it J. Math. Phys.}
{\bf 44} 1984-1997


\bibitem{book05} Giachetta G, Mangiarotti L and Sardanashvily G 2005 {\it
Geometric and Algebraic Topological Methods in Quantum Mechanics} (World
Scientific: Singapore)

\bibitem{guil} Guillemin V and Sternberg S 1984 {\it Symplectic Techniques
in Physics} (Cambr. Univ. Press: Cambridge)


\bibitem{karas} Karasev M and Maslov V 1993 Nonlinear Poisson brackets.
Geometry and quantization {\it Translations of AMS} {\bf 119} (AMS:
Providence, RI)

\bibitem{laz} Lazutkin V 1993 {\it  KAM Theory and Semiclassical
Approximations to  Eigenfunctions} (Springer: Berlin)

\bibitem{libe} Libermann P and Marle C-M 1987 {\it Symplectic Geometry and
Analytical Mechanics} (D.Reidel Publishing Company: Dordrecht)

\bibitem{magri} Magri F 1978 A simple model of the integrable Hamiltonian
equations {\it J. Math. Phys.} {\bf 19} 176-178

\bibitem{meig} Meigniez G 2002 Submersion, fibration and bundles {\it
Trans. Amer. Math. Soc.} {\bf 354} 3771-3787


\bibitem{mishc} Mishchenko A and Fomenko A 1978 Generalized Liouville
method of integration of Hamiltonian systems {\it Funct. Anal. Appl.}
{\bf 12} 113-121


\bibitem{nekh} Nekhoroshev N 1972 Action-angle variables and their
generalization {\it Trans. Mos. Math. Soc.} {\bf 26} 180-198

\bibitem{nekh94} Nekhoroshev N 1994 The
Poincar\'e-Lyapunov-Liouville-Arnol'd theorem {\it Funct. Anal. Appl.}
{\bf 28} (1994) 128-129

\bibitem{sard98} Sardanashvily G 1998 Hamiltonian time-dependent
mechanics {\it J. Math. Phys.} {\bf 39} 2714-2729

\bibitem{smirn} Smirnov R 1996 On the master symmetries related to a
certain classes of integrable hamiltonian systems {\it J. Phys. A} {\bf
29} 8133-8138

\bibitem{spar} Sparano G and Vilasi G 2000 Noncommutative integrability
and recursion operators {\it J. Geom. Phys.} {\bf 36} 270-284

\bibitem{stef} Stefan A 1974 Accesible sets, orbits and foliations with
singularities {\it Proc. London Math. Soc.} {\bf 29} 699-713

\bibitem{susm} Sussmann H 1973 Orbits of families of vector fields and
integrability of distributions {\it Trans. Amer. Math. Soc.} {\bf 180}
171-188

\bibitem{vin} Vinogradov A and Kupershmidt B 1977 The structures
of Hamiltonian mechanics {\it Russian Math. Surveys} {\bf 32} (4)
177-244

\bibitem{zung} Zung Nguyen Tien 2005 Torus actions and integrable
systems In vol. {\it Topological Methods in Theory of Integrable
Systems} Eds. A. Bolsinov, A.Fomenko and A.Oshemkov (Cambr. Sci.
Publ.: Cambridge); {\it Preprint} math.DS/0407455


\end{thebibliography}
\end{document}